\documentclass{article}


\newcounter{parentnumber}

\usepackage{amsfonts}
\usepackage{amssymb}
\usepackage{textcomp}
\usepackage{latexsym,amsmath}
\usepackage{graphics}

\renewcommand{\phi}{\varphi}
\newcommand{\be}{\begin{equation}}
\newcommand{\ee}{\end{equation}}
\newcommand{\ba}{\begin{eqnarray}}
\newcommand{\ea}{\end{eqnarray}}
\newcommand{\ban}{\begin{eqnarray*}}
\newcommand{\ean}{\end{eqnarray*}}

\newcommand{\nul}{{\bf0}}
\newcommand{\rd}{{\mathbb R}^d}
\newcommand{\zd}{{\mathbb Z}^d}
\newcommand{\td}{{\mathbb T}^d}
\renewcommand{\r}{{\mathbb R}}
\newcommand{\z} {{\mathbb Z}}
\newcommand{\cn} {{\mathbb C}}
\newcommand{\n} {{\mathbb N}}

\newcommand{\ddd}{,\dots,}
\renewcommand{\lll}{\left(}
\newcommand{\rrr}{\right)}
\newcommand{\h}{\widehat}
\newcommand{\w}{\widetilde}

\def\N{{{\Bbb N}}}
\def\Z{{{\Bbb Z}}}

\def\R{{\Bbb R}}

\def\vp{{\varphi}}

\def\s{{\sigma}}

\def\C{{\Bbb C}}

\def\){\right)}
\def\({\left(}

\def\supp{\operatorname{supp}}
\def\sinc{\operatorname{sinc}}
\def\mes{\operatorname{mes}}

\title{Approximation by multivariate \\ Kantorovich--Kotelnikov operators
\thanks{This research is supported by Volkswagen Foundation; the first author is also
supported by the project AFFMA that has received funding from the European Union's Horizon 2020 research and innovation
programme under the Marie Sklodowska-Curie grant agreement No 704030; the second  author is also
supported by grants from RFBR \# 15-01-05796-a,  St. Petersburg State University \#~9.38.198.2015.}}
\author{
Yu. Kolomoitsev$^{1, 2}$ and M. Skopina$^{3}$
}
\date{\small $^{1}$Universit\"at zu L\"ubeck,
Institut f\"ur Mathematik, L\"ubeck, Germany \\
\small $^{2}$Institute of Applied Mathematics and Mechanics of NAS of Ukraine, Slov'yans'k, Ukraine  \\
 $^{3}$St. Petersburg State University, Russia \\
kolomoitsev@math.uni-luebeck.de, skopina@ms1167.spb.edu}

\begin{document}
\maketitle

\begin{abstract}
%
Approximation properties of  multivariate Kantorovich-Kotelnikov type  operators generated by different band-limited functions
 are studied. In particular, a wide class of functions with  discontinuous Fourier transform is considered. The $L_p$-rate of convergence for these operators is given in terms of the classical moduli of smoothness.
Several examples of the Kantorovich-Kotelnikov  operators generated by the $\sinc$-function and its linear combinations are provided.

\end{abstract}

\bigskip


\medskip



\newtheorem{theo}{Theorem}
\newtheorem{lem}[theo]{Lemma}
\newtheorem {prop} [theo] {Proposition}
\newtheorem {coro} [theo] {Corollary}
\newtheorem {defi} [theo] {Definition}
\newtheorem {rem} [theo] {Remark}
\newtheorem {ex} [theo] {Example}

\newtheorem{theorempart}{Theorem}[theo]

\newtheorem{lemmapart}{Lemma}[theo]

\newtheorem{proppart}{Proposition}[theo]

\newcommand{\tocsecindent}{\hspace{0mm}}

\section{Introduction}

The Kantorovich--Kotelnikov operator is an operator of the form
\begin{equation}\label{KK}
  K_w(f,\vp ;x)=\sum_{k\in \Z}\(w\int_{\frac kw}^{\frac{k+1}w}f(u)du\)\vp(w x-k),\quad x\in \R,\quad w>0,
\end{equation}
where $f \,:\, \R\to \C$ is a locally integrable function and $\vp$ is an appropriate kernel satisfying certain "good" properties, as a rule $\vp$ is a  band-limited function or a function with a compact support, e.g., $B$-spline.
The operator $K_w$ was introduced in~\cite{BBSV}, although in other forms, it was known previously, see, e.g.,~\cite{Jia1, Jia2, LJC}. During the last years, in view of some important applications, this operator has drawn attention by many mathematicians and  has been especially actively studied~\cite{BM, CS, CV2, CV0, FCMV, Jia2, KM, KKS, KS, KS1, OT15, VZ1, VZ2}.

The operator $K_w$ has several advantages  over the generalized sampling operators 
\begin{equation}\label{SS}
  S_w (f,\vp ;x)=\sum_{k\in \Z}f\(\frac kw\)\vp(w x-k),\quad x\in \R,\quad w>0.
\end{equation}
First of all, using the averages
$w\int_{\frac kw}^{\frac{k+1}w}f(u)du$ instead of the sampled values $f(k/w)$ allows to deal with discontinues signals and reduce the so-called time-jitter errors. Note that the latter property is very useful in the theory of Signal and Image Processing. Moreover, unlike to the generalized sampling operators $S_w$,  the operator~\eqref{KK} is continuous in $L_p(\R)$ and, therefore,  provides better approximation order than $S_w$ in most cases.

%
%

In the literature, there are several generalizations and refinements of the Kantorovich--Kotelnikov operator $K_w$ (see, e.g.,~\cite{BDR, v58, Jia2,  KKS, KS, KS1, LJC, OT15, VZ1, VZ2}).  In this paper, we study approximation properties of the following multivariate analogue of~\eqref{KK}
\begin{equation}\label{gKK}
  Q_j (f,\vp,\w\vp; x)=\sum_{k\in \Z^d}\(m^j \int_{\R^d}f(u)\w\vp (M^ju+k) du\)\vp(M^j x+k),\quad x\in \R^d,\quad j\in \Z,
\end{equation}
where $M$ is a dilation matrix, $m=|{\rm det}\,M|$, and $\w\vp$ and $\vp$ are appropriate functions. Note that if $d=1$ and $\w\vp(x)=\chi_{[0,1]}(x)$ (the characteristic function of $[0,1]$), then~\eqref{gKK} represents the standard Kantorovich--Kotelnikov operator $K_{m^j}$.

Convergence and approximation properties of the operator~\eqref{gKK} have been actively studied by many authors (see~\cite{BDR, BM, CS, CV2, CV0, FCMV, v58, Jia1, Jia2, KKS, KS, KS1, LJC, OT15, VZ1, VZ2}). The most general results on estimates of the error of approximation by $Q_j$ have been obtained in the case of compactly supported $\vp$ and $\w\vp$. In particular in~\cite{Jia2} (see also~\cite{Jia1}) it was proved the following result: \emph{ if $\vp$ and $\w\vp$ are compactly supported, $\vp \in L_p(\R^d)$, $\w\vp\in L_q(\R^d)$, $1/p+1/q=1$, $M$ is an isotropic dilation matrix, and $Q_0$ reproduces polynomials of degree $n-1$, then for any $f\in L_p(\R)$,
$1\le p\le \infty$, and $n\in \N$ we have
\begin{equation}\label{Qvved}
  \Vert f-Q_j(f,\vp,\w\vp)\Vert_{L_p(\R^d)}\le C({p,d,n,\vp,\w\vp})\omega_n(f,\|M^{-j}\|)_p,
\end{equation}}
\emph{where $\omega_n(f,h)_p$ is the modulus of smoothness of order $n$.}

Concerning band-limited functions $\vp$, it turns out that approximation properties of the operators $Q_j$ have been studied mainly in the case where $\w\vp$ is some characteristic function (see, e.g.,~\cite{BBSV, CS, CV2, CV0, FCMV, KKS, OT15, VZ1, VZ2}) and, unlike to compactly supported functions $\vp$, there are several limitations and drawbacks in the available results.
First of all,  the methods previously used  are essentially restricted to the case of integrable functions $\vp$, which do not allow to consider the functions of type $\sinc(x)=(\sin \pi x)/(\pi x)$. Secondly, the conditions imposed on the kernel $\vp$  cannot provide high rate of convergence of $Q_j(f)$ even for sufficiently smooth functions $f$. At that, the corresponding estimates  are given in the terms of Lipschitz classes.
For example, it follows from~\cite{CV2} that \emph{for any $f\in L_p(\R)\cap {\rm Lip}(\nu)$, $1\le p\le \infty$, $0<\nu\le 1$, we have
\begin{equation}\label{intrImp}
  \bigg\Vert f-\frac12\sum_{k\in \Z}\(w\int_{\frac kw}^{\frac{k+1}w}f(u)du\)\sinc^2\(\frac{wx-k}{2}\)\bigg\Vert_{L_p(\R)}=\mathcal{O}(w^{-\nu}),\quad w\to +\infty.
\end{equation}
}

In this paper, we improve the mentioned drawbacks and study the rate of convergence of the operator~\eqref{gKK} for a wide class of band-limited functions $\vp$ including non-integrable  ones. In particular, estimation~(\ref{Qvved}) is proved for a large class of functions $\w\phi$ including both compactly supported and band-limited functions, provided that
$D^{\beta}(1-\h\phi\overline{\h{\w\phi}})(\textbf{0}) = 0$ for all $\beta\in\zd_+$, $\|\beta\|_{\ell_1}<n$ (see Theorems~\ref{th2} and~\ref{th2}\,$'$).
In the partial cases, this gives an answer to the question posed in~\cite{BBSV} about approximation properties of the following sampling series (see Section~6):
\begin{equation*}
 \sum_{k\in \Z}\(w\int_{\frac kw}^{\frac{k+1}w}f(u)du\)\sinc(w x-k),\quad x\in \R,\quad w>0.
\end{equation*}




The paper is organized as follows: in
Section~2 we introduce notation and give  some basic facts.
In Section~3 we consider approximation properties of some generalized sampling operators of type $Q_j$. These operators are defined similarly to the operators $Q_j$ but with appropriate tempered distribution in place of the function $\w\vp$ that, in particular, allows to include in the consideration  operators of type $S_w$.
The $L_p$-rate of convergence of such generalized sampling operators is given in terms of Fourier transform of $f$ and has several drawbacks, which we improve in the next sections. Section~4 is devoted to auxiliary results. In Section~5 we prove two main results  that provide estimates of the error of approximation by the operator $Q_j$ in terms of the classical moduli of smoothness. The results of this section can be considered as counterparts of the corresponding results from Section~3. Finally,  in Section~6 we consider some special cases  and provide a number of examples.   

\section{Notation and basic facts}
\label{notation}

$\n$ is the set of positive integers,
    $\r$  is the set of real numbers,
    $\cn$ is the set of complex numbers.
    $\rd$ is the
    $d$-dimensional Euclidean space,  $x = (x_1\ddd x_d)$ and $y =
    (y_1\ddd y_d)$ are its elements (vectors),
    $(x, y)~=~x_1y_1+~\dots~+x_dy_d$,
    $|x| = \sqrt {(x, x)}$, ${\bf0}=(0\ddd 0)\in\rd$;
      $B_r=\{x\in\rd:\ |x|\le r\}$, $\td=[-\frac 12,\frac 12]^d$;
    $\zd$ is the integer lattice
    in $\rd$, $\z_+^d:=\{x\in\zd:~x\geq~{\bf0}\}.$
    If $\alpha,\beta\in\zd_+$, $a,b\in\rd$, we set
    $[\alpha]=\sum\limits_{j=1}^d \alpha_j$,
    $\alpha!=\prod\limits_{j=1}^d(\alpha_j!),$
    $$\binom{\beta}{\alpha}=\frac{\beta!}{\alpha!(\beta-\alpha)!},\quad
    D^{\alpha}f=\frac{\partial^{[\alpha]} f}{\partial x^{\alpha}}=\frac{\partial^{[\alpha]} f}{\partial^{\alpha_1}x_1\dots
    \partial^{\alpha_d}x_d},$$
    $\delta_{ab}$~is the Kronecker delta.

A real $d\times d$ matrix $M$ whose
eigenvalues are bigger than 1 in modulus is called a  dilation matrix.
Throughout the paper we
consider that such a matrix $M$ is fixed and  $m=|{\rm det}\,M|$,
$M^*$ denotes the conjugate matrix to $M$.
Since the spectrum of the operator $M^{-1}$ is
located in  $B_r$,
where $r=r(M^{-1}):=\lim_{j\to+\infty}\|M^{-j}\|^{1/j}$
is the spectral radius of $M^{-1}$, and there exists at least
one point of the spectrum on the boundary of $B_r$, we have
\be
	\|M^{-j}\|\le {C_{M,\vartheta}}\,\vartheta^{-j},\quad j\ge0,
	\label{00++}
\ee
for every  positive number $\vartheta$  which is smaller in modulus
than any eigenvalue of  $M$.
In particular, we can take $\vartheta > 1$, then
	$$
	\lim_{j\to+\infty}\|M^{-j}\|=0.
	$$
%

$L_p$ denotes $L_p(\rd)$, $1\le p\le\infty$, with the usual norm $\Vert f\Vert_p=\Vert f\Vert_{L_p(\R^d)}$.
We say that  $\vp\in L_p^0$ if $\vp\in L_p$ and $\vp$ has a compact support.
We use $W_p^n$, $1\le p\le\infty$, $n\in\n$,
to denote  the Sobolev space on~$\rd$, i.e. the set of
functions whose derivatives up to order $n$ are in $L_p$,
with the usual Sobolev semi-norm given by
$$
\|f\|_{\dot W_p^n}=\sum_{[\nu]=n}\Vert D^\nu f\Vert_p.
$$

If $f, g$ are functions defined on $\rd$ and $f\overline g\in L_1$,
then  $\langle  f, g\rangle:=\int_{\rd}f\overline g$.
If $f\in L_1$,  then its Fourier transform is $\mathcal{F}f(\xi)=\widehat
f(\xi)=\int_{\rd} f(x)e^{-2\pi i
(x,\xi)}\,dx$.

If $\phi$ is a function defined on $\rd$, we set
$$
\phi_{jk}(x):=m^{j/2}\phi(M^jx+k),\quad j\in\z,\,\, k\in\rd.
$$

Denote by $\mathcal{S}$ the Schwartz class of functions defined on $\rd$.
    The dual space of $\mathcal{S}$ is $\mathcal{S}'$, i.e. $\mathcal{S}'$ is
    the space of tempered distributions.
    The basic facts from distribution theory
    can be found, e.g., in~\cite{Vladimirov-1}.
    Suppose $f\in \mathcal{S}$, $\phi \in \mathcal{S}'$, then
    $\langle \phi, f\rangle:= \overline{\langle f, \phi\rangle}:=\phi(f)$.
    If  $\phi\in \mathcal{S}',$  then $\h \phi$ denotes its  Fourier transform
    defined by $\langle \h f, \h \phi\rangle=\langle f, \phi\rangle$,
    $f\in \mathcal{S}$.
    If  $\phi\in \mathcal{S}'$, $j\in\z, k\in\zd$, then we define $\phi_{jk}$ by
        $
        \langle f, \phi_{jk}\rangle=
        \langle f_{-j,-M^{-j}k},\phi\rangle$ for all $f\in \mathcal{S}$.

Denote by $\mathcal{S}_N'$ the set of all tempered distribution
	whose Fourier transform $\h{\phi}$ is a function on $\rd$
	such that $|\h{\phi}(\xi)|\le C_{\phi} |\xi|^{N}$
	 for almost all $\xi\notin\td$, $N=N({\phi})\ge 0,$ and
	 $|\h{\phi}(\xi)|\le C'_{{\phi}}$
	 for almost all $\xi\in\td$.

Denote by ${\cal B}={\cal B}(\R^d)$  the class of functions $\phi$ given by
$$
\phi(x)=\int\limits_{\rd}\theta(\xi)e^{2\pi i(x,\xi)}\,d\xi,
$$
where $\theta$ is supported in a parallelepiped $\Pi:=[a_1, b_1]\times\dots\times[a_d, b_d]$ and such that	$\theta\big|_\Pi\in C^d(\Pi)$.

Let $1\le p \le \infty$. Denote by ${\cal L}_p$ the set
	$$
	{\cal L}_p:=
	\left\{
	\phi\in L_p\,:\, \|\phi\|_{{\cal L}_p}:=
	\bigg\|\sum_{k\in\zd} \left|\phi(\cdot+k)\right|\bigg\|_{L_p(\td)}<\infty
	\right\}.
	$$
	With the norm $\|\cdot\|_{{\cal L}_p}$, ${\cal L}_p$ is a Banach space.
	The simple properties are:
	${\cal L}_1=L_1,$
	$\|\phi\|_p\le \|\phi\|_{{\cal L}_p}$,
	$\|\phi\|_{{\cal L}_q}\le \|\phi\|_{{\cal L}_p}$
	for $1\le q \le p \le\infty.$ Therefore, ${\cal L}_p\subset L_p$
	and ${\cal L}_p\subset {\cal L}_q$ for $1\le q \le p \le\infty.$
	If $\phi\in L_p$ and compactly supported, then $\phi\in {\cal L}_p$ for $p\ge1.$
	If $\phi$ decays fast enough, i.e. there exist constants $C>0$
	and $\varepsilon>0$ such that
	$|\phi(x)|\le C( 1+|x|)^{-d-\varepsilon}$ for all $x\in\rd,$
	then $\phi\in {\cal L}_\infty$.
	
The modulus of smoothness $\omega_n(f,\cdot)_p$  of order $n\in \N$ for a function $f\in L_p$  is defined by
$$
\omega_n(f,h)_p=\sup_{|\delta|<h,\, \delta\in \R^d} \Vert \Delta_\delta^n f\Vert_{p},
$$
where
$$
\Delta_\delta^n f(x)=\sum_{\nu=0}^n (-1)^\nu \binom{n}{\nu} f(x+\delta \nu).
$$


\section{Preliminary results}
\label{scaleAppr}

Scaling operator
$\sum_{k\in\zd} \langle f, {\w\phi}_{jk}\rangle \phi_{jk}$
is a good tool of approximation for many  appropriate pairs of functions
$\phi, \w\phi$.
Let us consider the case, where $\w\phi$ is a
tempered distribution, e.g., the delta-function or a linear
combination of its derivatives. In this case the inner product
$\langle f, \w\phi_{jk}\rangle $ has meaning only for functions $f$ in $\mathcal{S}$.
To extend the class of functions $f$ one can  replace
$\langle f, {\w\phi}_{jk}\rangle $  by
$\langle \h f, \h{\w\phi_{jk}}\rangle$.
Approximation properties of such operators for certain classes
of distributions~$\w\phi$ and functions $\phi$ were studied, e.g., in~\cite{KKS} and~\cite{KS1}.

Repeating step-by-step the proof of Theorem~4 in~\cite{KS1}
and using Corollary~10 in~\cite{KKS}, it is easy to prove the following result.

\begin{prop}
\label{theoQjnew}
	Let $N\in\z_+,$    $\w\phi\in \mathcal{S}_N'$,  $\phi \in \mathcal{B} $.
	Suppose that  there exist $n\in\n$ and $\delta\in(0, 1/2)$ such that
	$\h\phi\h{\w\phi}$ is  boundedly differentiable up to order $n$ on
	$\{|\xi|<\delta\}$;  ${\rm supp\,} \h\phi\subset B_{1-\delta}$;
		$D^{\beta}(1-\h\phi\h{\w\phi})(0) = 0$ for $[\beta]<n$.
If	 $2\le p< \infty$, $1/p+1/q=1$,
$\gamma\in(N+d/p, N+d/p+\varepsilon)$, $\varepsilon>0,$ and
     \be
f\in L_p,\  \
       \h f\in L_q,\ \  \h f(\xi)=O(|\xi|^{-N-d-\varepsilon})\
\text{as }|\xi|\to\infty,
 \label{(c)}	
\ee
	then
\begin{equation*}
  \begin{split}
    	\bigg\|f-\sum_{k\in\zd} \langle \widehat{f}, \h{\w\phi_{jk}}\rangle \phi_{jk}\bigg\|_p^q\le
	C_1
	 \|M^{*-j}\|^{\gamma q}&\int\limits_{|M^{*-j}\xi|\ge\delta}
	|\xi|^{q\gamma}| \h f(\xi)|^q d\xi\\
      &+C_2 \|M^{*-j}\|^{nq}\int\limits_{|M^{*-j}\xi|<\delta}
	|\xi|^{qn}|\h f(\xi)|^q d\xi,
  \end{split}
\end{equation*}
where $C_1$ and $C_2$ do not depend on $j$ and $f$.
\end{prop}	



Proposition~\ref{theoQjnew} does not provide approximation order of $\sum_{k\in\zd} \langle \widehat{f}, \h{\w\phi_{jk}}\rangle \phi_{jk}$ better than $\|M^{*-j}\|^{n}$ even
for very smooth functions $f$. This can be fixed under stronger restrictions on $\vp$ given in the following definition.

\begin{defi}
\label{d1}
 A tempered distribution  $\w\phi$ and a function $\phi$ is said to be
 {\em strictly compatible} if there exists $\delta\in(0,1/2)$ such that
 $\overline{\h\phi}(\xi)\h{\w\phi}(\xi)=1$
 a.e. on $\{|\xi|<\delta\}$ and $\h\phi(\xi)=0$ a.e.
 on $\{|l-\xi|<\delta\}$ for all
 $l\in\z^d\setminus\{0\}$.
\end{defi}


%

\begin{rem}
\emph{ It is well known that the shift-invariant space generated by a function $\varphi$ has approximation order $n$ if and only if
the Strang--Fix condition of order $n$ is satisfied for $\varphi$ (that is $D^\beta\widehat\varphi(l)=0$ whenever $l\in\zd$, $ l\ne\nul$, $[\beta]<n$).
This fact appeared in the literature in different situations many times (see, e.g.,~\cite{BDR},  \cite{DB-DV1}, \cite{v58}, and~\cite[Ch.~3]{NPS}).
The condition  $D^{\beta}(1-\h\phi\h{\w\phi})(0) = 0$, $[\beta]<n$, is also a natural requirement for providing approximation order $n$ of scaling operators. This assumption  often appears
(especially in wavelet theory) in other terms, in particular, in terms of polynomial reproducing property (see~\cite[Lemma~3.2]{Jia1}).
It is clear that to provide an infinitely large approximation order, these conditions should be satisfied for any~$n$. Clearly, the latter holds for strictly compatible functions $\vp$ and~$\w\vp$.
 }


\emph{Supposing that $\h{\w\phi}(\xi)=1$ a.e. on $\{|\xi|<\delta\}$, it is easy to see that the simplest example of $\varphi$ satisfying Definition 2 is the tensor product of the $\rm sinc$ functions.}
\end{rem}

\begin{prop} {\sc \cite[Theorem~11]{KKS}}
\label{theoQjcomp}
	Let $N\in\z_+,$    $\w\phi\in \mathcal{S}_N'$,  $\phi \in \mathcal{B} $,
		$\w\phi$ and $\phi$ be strictly compatible.

\noindent
If	 $2\le p< \infty$, $1/p+1/q=1$,
$\gamma\in(N+d/p, N+d/p+\varepsilon)$, $\varepsilon>0,$ and a function $f$ satisfies~(\ref{(c)}), then
	\be
	\bigg\|f-\sum_{k\in\zd} \langle \widehat{f}, \h{\w\phi_{jk}}\rangle \phi_{jk}\bigg\|_p^q\le C
	 \|M^{*-j}\|^{\gamma q}  \int\limits_{|M^{*-j}\xi|\ge\delta}
	|\xi|^{q\gamma}| \h f(\xi)|^q d\xi,
	 \label{8}
	\ee
where $C$ does not depend on $j$ and $f$.
\end{prop}

Note that Proposition~\ref{theoQjcomp} is an analog of the following Brown's result~\cite{Brown}:
$$
\left|f(x)- \sum_{k\in\,\z} f(-2^{-j}k)\,{\rm sinc}(2^jx+k)\right|\le
C \int\limits_{|\xi|>2^{j-1}}
 |\h f(\xi)|d\xi,\quad x\in\r,
$$
whenever the Fourier transform of  $f$ is summable on $\r$.

There are several drawbacks in  Propositions~\ref{theoQjnew} and~\ref{theoQjcomp}. First,  they are proved only in the case $p\ge 2$. Second, there are additional restrictions on the function $f$.
Even in the case $\w\phi\in L_q^0$, where $\sum_{k\in\zd}\langle{f}, {\w\phi_{jk}}\rangle \phi_{jk}$ has meaning for every  $f\in L_p$,  the error estimate is obtained only for functions $f$ satisfying~(\ref{(c)}) with $N=0$.
Third, the error estimate is  given in the terms of decreasing of Fourier transforms unlike to common estimates, which usually are given in the terms of moduli of smoothness.

Below, under more restrictive conditions on $\w\vp$, we obtain analogues of the above propositions for all $f\in L_p$, $1\le p\le\infty$, and give the estimates of the error of approximation in terms of the classical moduli of smoothness.



\section{Auxiliary results}

	The following auxiliary statements will be useful for us.


\begin{lem}
\label{prop2}
Let either $1< p< \infty$, $\phi\in \cal B$ or $1\le p\le \infty$, $\vp \in \mathcal{L}_p$, and $a=\{a_k\}_{k\in\zd}\in\ell_p$. Then
\be
\label{201}
\bigg\|\sum_{k\in\zd} a_k \phi_{0k}\bigg\|_p\le C_{\phi, p}\|a\|_{\ell_p}.
\ee
\end{lem}

{\bf Proof.} In the case $\phi\in \cal B$, the proof of~\eqref{201} follows from~\cite[Proposition 9]{KKS}. Concerning the case $\vp \in \mathcal{L}_p$ see~\cite[Theorem 2.1]{v58}.~$\Diamond$

%

\begin{lem}\label{prop1Lq}
Let $f\in L_p$. Suppose that either $1<p<\infty$, $\phi\in \mathcal{B}$ or $1\le p\le\infty$, $\phi\in  L_q^0$, $1/p+1/q=1$. Then
\begin{equation}
\label{eqK1Lq}
  \bigg(\sum_{k\in\zd} |\langle f,{{\phi}_{0k}}\rangle|^p\bigg)^\frac 1p\le C'_{\phi, p}\|f\|_p.
\end{equation}
\end{lem}

{\bf Proof.} In the case $\phi\in \mathcal{B}$, the proof of~\eqref{eqK1Lq} see in~\cite[Proposition~6]{KKS}. The case $\phi\in  L_q^0$ follows from~\cite[Lemmas~4 and 5]{KS}.~$\Diamond$

%
\medskip

Combining the above lemmas, we obtain the following statement.

\begin{lem}\label{lemBound}
  Let $f\in L_p$, $1<p<\infty$, and $j\in \N$. Suppose $\vp\in \mathcal{B}$ and $\widetilde{\vp}\in \mathcal{B}\cup  \mathcal{L}_q$, $1/p+1/q=1$. Then
  \begin{equation}\label{eq0}
    \bigg\|\sum_{k\in\zd} \langle f,\w\phi_{jk}\rangle \phi_{jk}\bigg\|_p\le C\Vert f\Vert_p,
  \end{equation}
  where $C$ does not depend on $f$ and $j$.
\end{lem}

{\bf Proof.} If $\w\vp\in \mathcal{B}$, then the proof of~\eqref{eq0} directly follows from Lemmas~\ref{prop2} and~\ref{prop1Lq}.
In the case $\w\phi\in  \mathcal{L}_q$, we find  $g\in L_q$, $\|g\|_q\le1$, such that
\begin{equation}\label{eq+eq1}
     \bigg\|\sum_{k\in\zd} \big\langle f,\w\phi_{jk}\big\rangle \phi_{jk}\bigg\|_p=   \bigg|\bigg\langle\sum_{k\in\zd} \langle f,\w\phi_{jk}\rangle \phi_{jk}, g\bigg\rangle\bigg|=\bigg|\bigg  \langle f, \sum_{k\in\zd}  \langle\phi_{jk}, g\rangle \w\phi_{jk}\bigg\rangle\bigg|.
\end{equation}
It follows from the H\"older inequality and  Lemmas~\ref{prop2} and~\ref{prop1Lq} that
\begin{equation}\label{eq+eq2}
  \begin{split}
      \bigg|\bigg  \langle f, \sum_{k\in\zd}  \langle\phi_{jk}, g\rangle \w\phi_{jk}\bigg\rangle\bigg|&\le \|f\|_p\Big\|  \sum_{k\in\zd}  \langle\phi_{jk}, g\rangle \w\phi_{jk}  \Big\|_q\\
      &\le C_{\w\phi, q}\|f\|_p\bigg( \sum_{k\in\zd}|  \langle\phi_{jk}, g\rangle|^q\bigg)^{1/q}
       \le C_{\w\phi, q}{ C'}_{\phi, q}\|f\|_p.
   \end{split}
\end{equation}
Thus, combining \eqref{eq+eq1} and~\eqref{eq+eq2}, we obtain~(\ref{eq0}).~$\Diamond$

\bigskip

Similarly one can prove the following generalization of Lemma~\ref{lemBound} to the limiting cases $p=1,\infty$.

\bigskip

\noindent\textbf{Lemma~\ref{lemBound}$\,'$}
 \emph{ Let $f\in L_p$, $p=1$ or $p=\infty$, $j\in \N$. Suppose that}
 \smallskip

\noindent (i) {\it $\vp\in L_1$ and $\w\vp\in \mathcal{L}_\infty$ in the case $p=1$;}

\noindent (ii) {\it $\vp\in \mathcal{L}_\infty$ and $\w\vp\in L_1$ in the case $p=\infty$.}

\smallskip

\noindent\emph{Then
  \begin{equation*}
    \bigg\|\sum_{k\in\zd} \langle f,\w\phi_{jk}\rangle \phi_{jk}\bigg\|_p\le C\Vert f\Vert_p,
  \end{equation*}
where $C$ does not depend on $f$ and $j$.}

\bigskip

%
%


%

\begin{lem}\label{corE}
Let $N\in \Z_+$, $\w\phi\in \mathcal{S}_N'$, $\phi\in { \mathcal{ B}}\cup { \mathcal{ L}}_2$, and
$\w\phi$ and $\phi$ be strictly compatible.
If a function $f\in L_2$ is such that
its Fourier transform is supported in $\{\xi:\ |M^{*-j}\xi|<\delta\}$,
where $\delta$ is from Definition~\ref{d1}, then
\be
f=\sum_{k\in\zd} \langle \widehat{f}, \h{\w\phi_{jk}}\rangle \phi_{jk}\quad a.e.
\label{203}
\ee
\end{lem}


{\bf Proof.} In the case $\vp\in \mathcal{B}$, the proof of the lemma follows from~(\ref{8}).
In the case $\vp\in \mathcal{L}_2$, it follows from~\cite[ eq. (4.15)]{KS1} that~(\ref{8}) holds and, therefore, one has~\eqref{203}.~$\Diamond$

\bigskip

Let $\mathcal{B}_p^\s$, $1\le p\le\infty$, $\s>0$, denote the set of all entire functions $f$ in $\C^d$ which are of (radial) exponential type $\s$ and being restricted to $\R^d$ belong to $L_p$.


\begin{lem}\label{lem1} {\sc \cite[Theorem~3]{Wil}}
  Let $n\in \N$, $1\le p\le\infty$, $0<h<2\pi/\s$, $|\xi|=1$, and $P_\s\in \mathcal{B}_p^\s$. Then
  $$
  \Vert D^{n,\xi} P_\s\Vert_p\le \(\frac{\s}{2\sin (\s h/2)}\)^n \Vert \Delta_{h\xi}^n P_\s\Vert_p,
  $$
  where
  $$
  D^{n,\xi} f(x)=\mathcal{F}^{-1}((i\cdot,\xi)^n \widehat{f}(\cdot))(x).
  $$
\end{lem}

\begin{coro}\label{coroNS}
In terms of Lemma~\ref{lem1}, we have
\begin{equation}\label{NS++}
   \Vert P_\s\Vert_{\dot W_p^n}\le  C\s^n\omega_n(P_\s,\s^{-1})_p,
\end{equation}
where $C$ does not depend on $P_\s$.
\end{coro}

{\bf Proof.}
By Lemma~\ref{lem1}, for any $1\le j\le d$, we have that
\begin{equation*}\label{NS}
  \bigg\Vert \frac{\partial^n}{\partial x_j^n} P_\s\bigg\Vert_{p}\le C\s^n\Vert \Delta_{\s^{-1} e_j}^n P_\s\Vert_p\le   C\s^n\omega_n(P_\s,\s^{-1})_p,
\end{equation*}
which obviously implies~\eqref{NS++}.~$\Diamond$

\bigskip

We need several basic properties of the modulus of smoothness (see, e.g.,~\cite[Ch.~4]{Nik}).

\begin{lem}\label{lemmod}
  Let $f,g\in L_p$, $1\le p\le \infty$, and $n\in \N$. Then for any $\delta>0$, we have

\medskip

  \noindent {\rm (i)} $\omega_n(f+g,\delta)_p\le \omega_n(f,\delta)_p+\omega_n(g,\delta)_p$;

\medskip

  \noindent {\rm (ii)} $\omega_n(f,\delta)_p\le 2^n\Vert f\Vert_p$;

\medskip

  \noindent {\rm (iii)} $\omega_n(f,\lambda\delta)_p\le (1+\lambda)^n\omega_n(f,\delta)_p,\, \lambda>0$.

\end{lem}

Let us also recall the  Jackson-type theorem in $L_p$ (see, e.g.,~\cite[Theorem 5.2.1 (7)]{Nik} or~\cite[5.3.2]{Timan}).

\begin{lem}\label{lemJ}
  Let $f\in L_p$, $1\le p\le \infty$, $n\in \N$, and $\s>0$. Then there exists
$P_\s\in \mathcal{B}_p^\s$ such that
\begin{equation}\label{J}
  \Vert f-P_\s\Vert_p\le C\omega_n(f,1/\s)_p,
\end{equation}
where $C$ is a constant independent of $f$ and $P_\s$.

\end{lem}

\section{Main results}

Our main results are based on the following lemma.

\begin{lem}\label{lemJack}
Let $1\le p\le \infty$, $n\in \N$, $\phi\in ({\mathcal{ B}}\cup \mathcal{L}_2)\cap L_p$,
and $\w\phi\in L_q$, $1/p+1/q=1$. Suppose that the functions   $\vp$ and $\w\vp$ are strictly compatible
and there exists a constant $c=c(n,p,d,\vp,\w\vp)$ such that
  \begin{equation}\label{eq0++}
    \bigg\|\sum_{k\in\zd} \langle f,\w\phi_{jk}\rangle \phi_{jk}\bigg\|_p\le c\Vert f\Vert_p.
  \end{equation}
 Then	
\begin{equation}\label{eqJ1}
    \bigg\|f-\sum_{k\in\zd} \langle f,\w\phi_{jk}\rangle \phi_{jk}\bigg\|_p\le C \omega_n\(f,\Vert M^{-j}\Vert\)_p,
  \end{equation}
where $C$ does not depend on $f$ and $j$.
%
\end{lem}

{\bf Proof.}
Let $g\in L_p\cap L_2$ be such that
\begin{equation}\label{eq*}
  \Vert f-g\Vert_p\le \omega_n\(f,\Vert M^{*-j}\Vert\)_p.
\end{equation}
Using~\eqref{eq0++} and~\eqref{eq*}, we have
\begin{equation}\label{eqJ2}
  \begin{split}
    \bigg\|f-\sum_{k\in\zd} \langle f,\w\phi_{jk}\rangle \phi_{jk}\bigg\|_p
    &\le \Vert f-g\Vert_p+\bigg\|g-\sum_{k\in\zd} \langle g,\w\phi_{jk}\rangle \phi_{jk}\bigg\|_p+\bigg\|\sum_{k\in\zd} \langle g-f,\w\phi_{jk}\rangle \phi_{jk}\bigg\|_p\\
    &\le (1+c)\omega_n\(f,\Vert M^{*-j}\Vert\)_p+\bigg\|g-\sum_{k\in\zd} \langle g,\w\phi_{jk}\rangle \phi_{jk}\bigg\|_p.
  \end{split}
\end{equation}

By Lemma~\ref{lemJ}, for any $n\in \N$ and $g\in L_p\cap L_2$ there exists a function $J_n(g)\,:\, \R^d\to \mathbb{C}$ such that $\supp\widehat{J_n(g)}\subset \{|\xi|<{\delta}{\Vert M^{*-j}\Vert^{-1}}\}$ and
\begin{equation}\label{eqJ3}
  \Vert g-J_n(g)\Vert_p \le C_1\omega_n\(g,\delta^{-1}\Vert M^{*-j}\Vert\)_p,
\end{equation}
where $\delta$ is from Definition~\ref{d1} and $C_1$ does not depend on $g$ and $j$.

By Lemma~\ref{corE}, taking into account that $ \langle \widehat{J_n(g)},\widehat{\w\phi_{jk}}\rangle= \langle J_n(g),\w\phi_{jk}\rangle$ and $\{|\xi|<{\delta}{\Vert M^{*-j}\Vert^{-1}}\}\subset \{|M^{*-j}\xi|<\delta\}$, we have
\begin{equation}\label{eqJ4}
  \sum_{k\in\zd} \langle J_n(g),\w\phi_{jk}\rangle \phi_{jk}=J_n(g).
\end{equation}
Thus, using~\eqref{eqJ3}, \eqref{eqJ4}, and~\eqref{eq0++}, we derive
\begin{equation}\label{eqJ5}
  \begin{split}
    \bigg\|g-\sum_{k\in\zd} \langle g,\w\phi_{jk}\rangle \phi_{jk}\bigg\|_p &\le \Vert g-J_n(g)\Vert_p+\bigg\|\sum_{k\in\zd} \langle g-J_n(g),\w\phi_{jk}\rangle \phi_{jk}\bigg\|_p\\
    &\le (1+c)\Vert g-J_n(g)\Vert_p\le C_2\omega_n\(g,\delta^{-1}\Vert M^{*-j}\Vert\)_p.
  \end{split}
\end{equation}
Next, using Lemma~\ref{lemmod} and~\eqref{eq*}, we get
\begin{equation}\label{eqJ6}
  \begin{split}
    \omega_n\(g,\delta^{-1}\Vert M^{*-j}\Vert\)_p&\le (1+\delta^{-1})^n\omega_n\(g,\Vert M^{*-j}\Vert\)_p\\
    &\le C_3\(\Vert f-g\Vert +\omega_n\(f,\Vert M^{*-j}\Vert\)_p\)\le 2C_3\omega_n\(f,\Vert M^{*-j}\Vert\)_p\\
    &=C_4\omega_n\(f,\Vert M^{-j}\Vert\)_p.
  \end{split}
\end{equation}

Finally, combining~\eqref{eqJ2}, \eqref{eqJ5}, and~\eqref{eqJ6}, we obtain~\eqref{eqJ1}.~$\Diamond$

\bigskip

The following  statement  is a multivariate analogue of Theorem~7 in~\cite{Sk1}.
It is also a counterpart of Proposition~\ref{theoQjcomp} in some sense.

\begin{theo}\label{thJack}
Let $f\in L_p$, $1<p<\infty$, and $n\in \N$. Suppose $\vp\in \mathcal{B}$,
$\w\phi\in \mathcal{B}\cup \mathcal{L}_q$, $1/p+1/q=1$, and $\phi$ and $\w\phi$ are strictly compatible.
Then	
\begin{equation}
\label{eqJ1++}
    \bigg\|f-\sum_{k\in\zd} \langle f,\w\phi_{jk}\rangle \phi_{jk}\bigg\|_p\le C \omega_n\(f,\Vert M^{-j}\Vert\)_p,
  \end{equation}
where $C$ does not depend on $f$ and $j$.
%
\end{theo}

{\bf Proof.}
The proof follows from Lemma~\ref{lemJack} and Lemma~\ref{lemBound}.
$\Diamond$

\bigskip

\begin{rem}
\emph{By \eqref{00++}, it is easy to see that  in~\eqref{eqJ1++} and in further similar estimates, the modulus $\omega_n\(f,\Vert M^{-j}\Vert\)_p$ can be replaced by $\omega_n\(f,\vartheta^{-j}\)_p$, where $\vartheta$  is a positive number smaller in modulus
than any eigenvalue of $M$. Note that $\| M^{-1}\|$ may be essentially bigger than $\vartheta^{-1}$.}
\end{rem}

In the following result, we give an analogue of Theorem~\ref{thJack} for the limiting cases  $p=1$ and $p=\infty$.

\medskip


\noindent \textbf{Theorem~\ref{thJack}\,$'$}
 {\it Let $f\in L_p$, $p=1$ or $p=\infty$, and $n\in \N$. Suppose the functions $\phi$ and $\w\phi$ are strictly compatible and}

\medskip

\noindent (i) {\it  $\vp\in \mathcal{B}\cap L_1$ and $\w\vp\in \mathcal{L}_\infty$ in the case $p=1$;}

\medskip

\noindent (ii) {\it $\vp\in \mathcal{L}_\infty$ and $\w\vp\in L_1$ in the case $p=\infty$.

\medskip

\noindent  Then
\begin{equation*}
    \bigg\|f-\sum_{k\in\zd} \langle f,\w\phi_{jk}\rangle \phi_{jk}\bigg\|_p\le C \omega_n\(f,\Vert M^{-j}\Vert\)_p,
  \end{equation*}
where $C$ does not depend on $f$ and $j$.}

\bigskip

\textbf{Proof.} The proof follows from Lemmas~\ref{lemJack} and~\ref{lemBound}$\,'$.~$\Diamond$

\bigskip

\begin{rem}
 It is easy to see that Theorem~\ref{thJack}\,$'$ is valid if we replace the condition $\vp\in \mathcal{B}\cap L_1$ by $\vp\in \mathcal{L}_2$.
\end{rem}

The next theorem is the main result of the paper. This result essentially extends the classes of functions $\vp$ and $\w\vp$ from  Theorem~\ref{thJack} and can be considered as a counterpart of Proposition~\ref{theoQjnew}.

\begin{theo}\label{th2}
Let $f\in L_p$, $1<p<\infty$, and $n\in\n$.  Suppose
 $\phi\in\cal B$, ${\rm supp\,} \h\phi\subset B_{1-\varepsilon}$ for some $\varepsilon\in (0,1)$, $\h\phi\in C^{n+d+1}(B_\delta)$ for some $\delta>0$;
 $\w\phi \in\mathcal{B}\cup \mathcal{L}_q$, $1/p+1/q=1$,  $\h{\w\phi}\in C^{n+d+1}(B_\delta)$  and  $D^{\beta}(1-\h\phi\overline{\h{\w\phi}})({\bf 0}) = 0$ for all $\beta\in\zd_+$, $[\beta]<n$.  Then
\begin{equation}\label{1}
  \bigg\|f-\sum\limits_{k\in\zd}
\langle f,\widetilde\phi_{jk}\rangle \phi_{jk}\bigg\|_p\le C \omega_n\(f,\|M^{-j}\|\)_p,
\end{equation}
where $C$ does not depend on $f$ and $j$.
\end{theo}

\textbf{Proof.}
First, let us prove that for any $f\in W_p^n$
\begin{equation}\label{2}
  \bigg\|f-\sum\limits_{k\in\zd}
\langle f,\widetilde\phi_{jk}\rangle \phi_{jk}\bigg\|_p\le C_{0}\|f\|_{\dot W_p^n}\|M^{-j}\|^n,
\end{equation}
where $C_{0}$ does not depend on $f$ and $j$.

 Choose $0<\delta'<\delta''$ such that $\h\phi(\xi)\ne0$ on $\{|\xi|\le\delta'\}$ and $\delta'\le\delta$. Set
$$
F(\xi)=
 \begin{cases}
	 \displaystyle\frac{1-\overline{\h\phi(\xi)}\h{\w\phi}(\xi)}{\overline{\h\phi(\xi)}}  &\mbox{if $|\xi|\le\delta'$,}
		 \\
	\displaystyle 0
	 &\mbox{if $|\xi|\ge\delta''$}
	 		\end{cases}
	$$
and extend this function  such that $F\in C^{n+d+1}(\rd)$. Define $\w\psi$ by
$\h{\w\psi}=F$. Obviously, the function $\w\psi$ is continuous and $\w\psi(x)=O(|x|^{-\gamma})$ as $|x|\to \infty$, where
$\gamma>n+d$.

Since $\w\psi\in \cal B$ and $\w\phi\in \mathcal{L}_q\cup \mathcal{B}$, by Lemma~\ref{lemBound}, we have
\begin{equation*}
  \begin{split}
     \bigg\Vert\sum_{k\in\zd}\langle f, \w\phi_{jk}+\w\psi_{jk}\rangle \vp_{jk}\bigg\Vert_p&\le
\bigg\Vert\sum_{k\in\zd}\langle f, \w\phi_{jk}\rangle \vp_{jk}\bigg\Vert_p+
\bigg\Vert\sum_{k\in\zd}\langle f, \w\psi_{jk}\rangle \vp_{jk}\bigg\Vert_p\le C_1\|f\|_p.
   \end{split}
\end{equation*}
Now, taking into account that  $\overline{\h\phi(\xi)}(\h{\w\phi}(\xi)+\h{\w\psi}(\xi))=1$ whenever $|\xi|\le\delta'$,
we obtain from Lemma~\ref{lemJack} that for every $f\in W_p^n$
$$
\bigg\|f-\sum\limits_{k\in\zd}
\langle f,\widetilde\phi_{ik}+\widetilde\psi_{jk}\rangle \phi_{jk}\bigg\|_p
\le
C_2\omega_n\(f, \|M^{-j}\|\)_p\le C_3\|f\|_{\dot W_p^n}\|M^{-j}\|^n.
$$
Thus, to prove~\eqref{2},
it remains to verify that for  $f\in W_p^n$
\begin{equation}\label{dop1}
  \bigg\|\sum\limits_{k\in\zd}\langle f,\widetilde\psi_{jk}\rangle \phi_{jk}\bigg\|_p
\le C_4\|f\|_{\dot W_p^n}\|M^{-j}\|^n.
\end{equation}

Let $k\in\zd$, $z\in [-1/2,1/2]^d-k$, $y=M^{-j}z$. Since $D^{\beta}\h{\w\psi}(0) = 0$ whenever $[\beta]<n$, we have
$$
\int\limits_{\rd}y^\alpha\widetilde\psi_{jk}(y)\,dy=0,\quad j\in\z,\quad
 \alpha\in\zd_+,\quad [\alpha]< n.
$$
Hence, due to Taylor's formula with the integral remainder,
\begin{equation*}
  \begin{split}
      |\langle f,\w\psi_{jk}\rangle |&=
    \bigg|\int\limits_{\rd}f(x)\overline{\w\psi_{jk}(x)}\,dx\bigg|
    \\
    &=\bigg|\int\limits_{\rd}\, \overline{\w\psi_{jk}(x)}\bigg(\sum\limits_{\nu=0}^{n-1}\frac{1}{\nu!}\lll (x_1-y_1)\partial_{1}
    + \dots + (x_d-y_d)\partial_{d} \rrr^\nu f(y)
    \\
    &\quad\quad\quad\quad+\int\limits_0^1\frac{(1-t)^{n-1}}{(n-1)!} \Big((x_1-y_1)\partial_{1} + \dots + (x_d-y_d)\partial_{d} \Big)^n f(y+t(x-y))\,dt\bigg)\,dx \bigg|
    \\
    &\le \int\limits_{\rd}\,dx\,|x-y|^n|\w\psi_{jk}(x)|\,\int\limits_0^1\sum\limits_{[\beta]=n}|D^{\beta}f(y+t(x-y))|\,dt.
   \end{split}
\end{equation*}
From this, using  H\"older's inequality and taking into account that		
$$
|x-y|^n\le\|M^{-j}\|^n|M^j x-z|^n,
$$
$$
|\w\psi_{jk}(x)|\le\frac{C_5 m^{j/2}}{{(1+|M^jx+k|)^\gamma}}\le \frac{C_6 m^{j/2}}{(1+|M^jx-z|)^\gamma},
$$
we obtain
\begin{equation}\label{dlin}
  \begin{split}
          |\langle f,\w\psi_{jk}\rangle |&\le
   C_6  m^{j/2}\|M^{-j}\|^n \int\limits_{\rd}\,dx \frac{|M^j x-z|^n}{(1+|M^jx-z|)^\gamma}
    \int\limits_0^1\sum\limits_{[\beta]=n}|D^{\beta}f(y+t(x-y))|\,dt
    \\
    &\le C_6  m^{j/2}\|M^{-j}\|^n
\bigg(\int\limits_{\rd}\frac{|M^j x-z|^{n}}{(1+|M^jx-z|)^\gamma}\,dx\bigg)^{1/q}
\\
    &\quad\quad\quad\quad\times\Bigg(\int\limits_{\rd}\,\frac{|M^j x-z|^{n} dx}{(1+|M^jx-z|)^\gamma}\bigg(\int\limits_0^1\,
\sum\limits_{[\beta]=n}|D^{\beta}f(y+t(x-y))|\,dt\bigg)^p\Bigg)^{\frac 1 p}
\\
	&=C_6  m^{\frac j2-\frac jq}\|M^{-j}\|^n \bigg(\int\limits_{\rd}\frac{|x-z|^{n} \,dx}{(1+|x-z|)^\gamma}\bigg)^{1/q}
\\
	&\quad\quad\quad\quad\times\Bigg(\int\limits_{\rd}\,\frac{|M^j x-z|^{n} dx}{(1+|M^jx-z|)^\gamma}\bigg(\int\limits_0^1\,
    \sum\limits_{[\beta]=n}|D^{\beta}f(y+t(x-y))|\,dt\bigg)^p\Bigg)^{\frac 1 p}
\\
    &\le C_7  m^{\frac j2-\frac jq}\|M^{-j}\|^n
    \Bigg(\int\limits_{\rd}\,\frac{|M^j(x-y)|^{n} dx}{(1+|M^j(x-y)|)^\gamma}\int\limits_0^1\,\sum\limits_{[\beta]=n}|D^{\beta}f(y+t(x-y))|^p\,dt\Bigg)^{\frac 1 p}
\\
	&=C_7 m^{\frac j2-\frac jq}\|M^{-j}\|^n
    \Bigg(\int\limits_{\rd}\,\frac{|M^ju|^n du}{(1+|M^ju|)^\gamma}\int\limits_0^1\,\sum\limits_{[\beta]=
	n}|D^{\beta}f(y+tu)|^p\,dt\Bigg)^{\frac 1 p}.
   \end{split}
\end{equation}
Next, it follows from~\eqref{dlin} and Lemma~\ref{prop2} that
\begin{equation*}
  \begin{split}
         \bigg\|\sum\limits_{k\in\,\zd}\langle
    f,\w\psi_{jk}\rangle\phi_{jk}\bigg\|_p^p
    &\le
    C_8 m^{p(\frac j2-\frac jp)}\sum_{k\in\zd}|\langle f,\w\psi_{jk}\rangle |^p
\\
    &=C_8m^{p(\frac j2-\frac jp)}\sum_{k\in\zd}\int\limits_{[-1/2,1/2]^d-k}\,dz|\langle f,\w\psi_{jk}\rangle |^p
\\
    &\le C_9\|M^{-j}\|^{pn} \int\limits_{\rd}dz
    \int\limits_{\rd}\,\frac{|M^ju|^n  du}{(1+|M^ju|)^\gamma}\int\limits_0^1\,\sum\limits_{[\beta]=n}|D^{\beta}f(M^{-j}z+tu)|^p\,dt
\\
    &=C_9 \|M^{-j}\|^{pn} \int\limits_{\rd}dz
    \int\limits_{\rd}\,\frac{|u|^n du}{(1+|u|)^\gamma}\int\limits_0^1\,\sum\limits_{[\beta]=n}|D^{\beta}f(z+tM^{-j}u)|^p\,dt
\\
	&=C_9 \|M^{-j}\|^{pn}
    \int\limits_{\rd}\,\frac{|u|^n du}{(1+|u|)^\gamma}\int\limits_0^1\,dt \int\limits_{\rd}\sum\limits_{[\beta]=n}|D^{\beta}f(z+tM^{-j}u)|^p\,dz
\\
    &=C_9 \|M^{-j}\|^{pn}
    \int\limits_{\rd}\,\frac{|u|^n du}{(1+|u|)^\gamma}\int\limits_0^1\,dt\int\limits_{\rd}\sum\limits_{[\beta]=n}|D^{\beta}f(z)|^p\,dz
\\
    &\le C_{10}\|M^{-j}\|^{pn}\|f\|_{\dot W^n_p}^p.
   \end{split}
\end{equation*}
This implies~\eqref{dop1} and, therefore,~\eqref{2}.

%
%

Now let us prove inequality~\eqref{1}.
By Lemma~\ref{lemJ}, there exists $P_\s\in \mathcal{B}_p^\s$, $\s=1/\|M^{-j}\|$,  such that
\begin{equation}\label{J}
  \Vert f-P_\s\Vert_p\le C_{11}\omega_n(f,1/\s)_p.
\end{equation}
Using Lemma~\ref{lemBound}, we obtain
\begin{equation}\label{1+}
  \begin{split}
    \bigg\|f-\sum\limits_{k\in\zd}
\langle f,\widetilde\phi_{jk}\rangle \phi_{jk}\bigg\|_p\le C_{12}\Vert f-P_\s\Vert_p+\bigg\|P_\s-\sum\limits_{k\in\zd}
\langle P_\s,\widetilde\phi_{jk}\rangle \phi_{jk}\bigg\|_p.
  \end{split}
\end{equation}
Next, by~\eqref{2} and Corollary~\ref{coroNS}, we derive
\begin{equation}\label{2+}
  \begin{split}
    \bigg\|P_\s-\sum\limits_{k\in\zd}
\langle P_\s,\widetilde\phi_{jk}\rangle \phi_{jk}\bigg\|_p\le C_0\|P_\s\|_{\dot W_p^n}\|M^{-j}\|^n\le C_{13}\omega_n\(P_\s,1/\s\)_p.
  \end{split}
\end{equation}
Finally, combining~\eqref{J}--\eqref{2+} and using Lemma~\ref{lemmod}, we get
\begin{equation*}
  \begin{split}
    \bigg\|f-\sum\limits_{k\in\zd}
\langle f,\widetilde\phi_{jk}\rangle \phi_{jk}\bigg\|_p&\le C_{14}\(\Vert f-P_\s\Vert_p+\omega_n(P_\s,1/\s)_p\)\\
&\le C_{15}\(\Vert f-P_\s\Vert_p+\omega_n(f,1/\s)_p\)\le C\omega_n\(f,\|M^{-j}\|\)_p,
  \end{split}
\end{equation*}
which implies~\eqref{1}.~$\Diamond$

\bigskip

In the following result, we give an analogue of Theorem~\ref{th2} for the limiting cases  $p=1$ and $p=\infty$.

\bigskip

\noindent \textbf{Theorem~\ref{th2}\,$'$}
 {\it Let $f\in L_p$, $p=1$ or $p=\infty$, and $n\in \N$. Suppose the functions $\vp$ and $\w\vp$ are such that $\h\phi, \h{\w\phi}\in C^{n+d+1}(B_\delta)$  for some $\delta>0$, $D^{\beta}(1-\h\phi\overline{\h{\w\phi}})({\bf 0}) = 0$ for all $\beta\in\zd_+$, $[\beta]<n$, $\vp\in \mathcal{B}$, \, ${\rm supp\,} \h\phi\subset B_{1-\varepsilon}$ for some $\varepsilon\in (0,1)$, and}

\medskip

\noindent (i) {\it $\vp\in L_1$ and  $\w\vp\in \mathcal{L}_\infty$ in the case  $p=1$;}

\medskip

\noindent (ii) {\it $\vp\in \mathcal{L}_\infty$ and $\w\vp\in L_1$ in the case $p=\infty$. }

\medskip

\noindent{\it  Then

\begin{equation*}
  \bigg\|f-\sum\limits_{k\in\zd}
\langle f,\widetilde\phi_{jk}\rangle \phi_{jk}\bigg\|_p\le C \omega_n\(f,\|M^{-j}\|\)_p,
\end{equation*}
where $C$ does not depend on $f$ and $j$.
}

%
%
%
%
%

\bigskip

\textbf{Proof.} The proof is similar to the proof of Theorem~\ref{th2}. We only note that one needs to use Lemma~\ref{lemBound}\,$'$ instead of Lemma~\ref{lemBound}. At that, in the case $p=\infty$, using the first inequality in~\eqref{dlin} and Lemma~\ref{prop2}, we get
\begin{equation*}
   \begin{split}
         \bigg\|\sum\limits_{k\in\,\zd}\langle
    f,\w\psi_{jk}\rangle\phi_{jk}\bigg\|_\infty
    &\le
    C_1 m^{\frac j2}\sup_{k\in\zd}|\langle f,\w\psi_{jk}\rangle |\\
    &\le
    C_2  m^{j}\|M^{-j}\|^n \Vert f\Vert_{\dot W_\infty^n} \sup_{k\in\zd} \int\limits_{\rd} \frac{|M^j x-z|^n}{(1+|M^jx-z|)^\gamma}\,dx\\
    &\le
    C_3  \|M^{-j}\|^n \Vert f\Vert_{\dot W_\infty^n}.
   \end{split}
\end{equation*}
Theorem~\ref{th2}\,$'$ is proved.~$\Diamond$

\bigskip

\begin{rem}
Note that any function $\w\vp\in L_p^0$, $1\le  p\le \infty$, satisfies assumptions $\w\vp\in \mathcal{L}_p\subset L_p$ and $\h{\w\phi}\in C^{n+d+1}(B_\delta)$, and hence can be used in Theorems~\ref{th2} and~\ref{th2}\,$'$.
\end{rem}


\section{Special Cases }

I. Let us start from the classical multivariate Kotelnikov type decomposition that can be obtain by using the function
$$
\sinc(x):=\prod_{\nu=1}^d \frac{\sin(\pi x_\nu)}{\pi x_\nu},\quad x\in \R^d.
$$

In what follows, we restrict ourselves by the case $1<p<\infty$.

\begin{prop}
\label{pr1 }
Let $f\in L_p$,  $1<p<\infty$, and let $U$ be a bounded measured subset of $\rd$. Then
\begin{equation}\label{101}
  \bigg\Vert f-\sum_{k\in\zd}\frac{m^j}{\mes U}\int\limits_{M^{-j}{U}}f(-M^{-j}k+t)\,dt\,\sinc(M^j\cdot+k)\bigg\Vert_p\le  C \omega_1(f,\|M^{-j}\|)_p,
\end{equation}
where the constant $C$ does not depend on $f$ and $j$.
If, in addition, $U$ is symmetric with respect to the origin, then in~\eqref{101} the modulus of continuity $\omega_1(f,\|M^{-j}\|)_p$ can be replaced by the second order modulus of smoothness $\omega_2(f,\|M^{-j}\|)_p$.
\end{prop}

\textbf{Proof.} We use Theorem~\ref{th2} for
$$
 \phi(x)=\sinc(x)\quad\text{and}\quad\w\vp (x)=\frac1{\mes U}\chi_U(x).
$$
Then, taking into account that $\h{\w\phi}(\nul)=1$ and
$$
\langle f,\widetilde\phi_{jk}\rangle=\frac{m^{j/2}}{\mes U}\int\limits_{M^{-j}U}f(-M^{-j}k+t)\,dt,
$$
we can verify that all assumptions of Theorem~\ref{th2} are satisfied  with $n=1$, which provides inequality~\eqref{101}.

Now, let $U$ be symmetric with respect to the origin. In this case,
$$
\frac{\partial}{\partial x_j}(1-\h\phi\overline{\h{\w\phi}})(\textbf{0})=-\frac{\partial \overline{\h{\w\phi}}}{\partial x_j}(\textbf{0})=\frac{2\pi i}{\mes U}\int_U x_j dx=0,\quad 1\le j\le d.
$$
Therefore, all assumptions of Theorem~\ref{th2}  are satisfied  with $n=2$. 
$\Diamond$

\medskip

\begin{rem}
Relation~(\ref{101}) gives a general answer to the question posed in~\cite{BBSV} concerning approximation properties of the sampling series given by
$$
\sum_{k\in\zd}\frac{m^j}{\mes U}\int\limits_{M^{-j}{U}}f(-M^{-j}k+t)\,dt\,\sinc(M^j\cdot+k)
$$
in the spaces $L_p$ for $1<p<\infty$.
\end{rem}

\begin{rem}
It follows from Theorem~\ref{th2}\,$'$ that Proposition~\ref{pr1 } is valid for all $f\in L_p$, $1\le p\le\infty$, if we replace $\sinc(x)$ by $\sinc^2(x)$ in~\eqref{101}.
In particular, this gives an improvement of estimate~\eqref{intrImp}. The same conclusion holds for all propositions presented below.
\end{rem}

%
%
%

\bigskip

II. Now let us show that using of an appropriate linear combination of the function $\sinc(x)$ rather than this function itself can provide better rates of the approximation by the corresponding sampling operator.

\begin{prop}
\label{pr2}
Let $f\in L_p$, $1<p<\infty$, $n\in \N$, and let $U$ be a bounded measured subset in $\rd$. Then there exists a finite set of numbers $\{a_l\}_{l\in \Z^d}\subset \C$ depending only on $d$, $n$, and $U$ such that
for
\begin{equation}\label{eqfunc}
  \phi(x)=\sum_l a_l \sinc(x+l)
\end{equation}
we have
\begin{equation}\label{eqcor1++}
  \bigg\Vert f-\sum_{k\in\zd} \frac{m^{j}}{\mes U} \int\limits_{M^{-j}U}f(-M^{-j}k+t)\,dt\,\phi(M^j\cdot+k)\bigg\Vert_p\le  C \omega_n(f,\|M^{-j}\|)_p,
\end{equation}
where the constant $C$ does not depend on $f$ and $j$.
\end{prop}

\textbf{Proof.}
Let $\w\vp (x)=\frac1{\mes U}\chi_U(x)$.
Find complex numbers $c_\alpha$, $\alpha\in\zd_+$, $[\alpha]< n$, satisfying
$$
c_\nul=1, \quad\sum\limits_{\nul\le\alpha\le\beta}\left({\beta\atop\alpha}\right)
\overline{D^{\beta-\alpha}\h{\w\phi}(\nul)}c_\alpha=0\quad \forall \beta\in\zd_+,\, \nul<[\beta]< n
$$
and set
\be
T(\xi)=\sum_{\nul\le[\alpha]\le n}c_\alpha \prod_{j=1}^dg_{\alpha_j}(\xi_j),
\label{102}
\ee
where $g_k$ is a trigonometric polynomial such that $\frac{d^l g_k}{dt^l}(0)=\delta_{kl}$ for all $l=0,\dots, k$.
It is not difficult to deduce explicit recursive formulas for finding such polynomials (see, e.g.,~\cite[Sec. 3.4]{KPS} or~\cite{Sk3}).
Obviously, $D^\alpha T(\nul)=c_\alpha$ and $D^\alpha(T\cdot\h\phi)(\nul)=D^\alpha(T\cdot\chi_{[-1/2,1/2]^d})(\nul)=c_\alpha$ for all $\alpha\in\zd_+$, $[\alpha]< n$.
If now $T(\xi)=\sum_l a_l e^{2\pi i(l,\xi)}$, then setting
$
  \phi(x)=\sum_l a_l \sinc(x+l),
$
we obtain that  $D^{\beta}(1-\h\phi\overline{\h{\w\phi}})(\textbf{0}) = 0$ for all $\beta\in\zd_+$, $[\beta]<n$. Thus, due to Theorem~\ref{th2}, we have inequality~\eqref{eqcor1++}.
$\Diamond$

\bigskip

Let us write explicit formulas for the function~\eqref{eqfunc} and for the polynomial $T$ given by~\eqref{102} in the cases $d=1, 2$ and $n=4$.

\medskip

\textbf{Example 1.} Let first $d=1$, $U=[-1/2, 1/2]$, and $M=2$. Then
$$
\h{\w\phi}(0)=1,\  {\h{\w\phi}}'(0)=0,\ {\h{\w\phi}}''(0)=-\frac{\pi^2}{3},\ {\h{\w\phi}}'''(0)=0,
$$
which yields
$
c_0=1,\ c_1=0,\ c_2=\frac{\pi^2}{3},\  c_3=0,
$
and
$$
T(\xi)=1+\frac{\pi^2}{3}g_2(\xi),
$$
where
$$
g_2(u)=-\frac1{8\pi^2}(2-5e^{2\pi i u}+4e^{4\pi i u}-e^{6\pi iu}).
$$
Hence
$$
\phi(x)=\frac{11}{12}\sinc(x)+\frac{5}{24}\sinc(x+1)-\frac{1}{6}\sinc(x+2)+\frac{1}{24}\sinc(x+3)
$$
and, by Proposition~\ref{pr2}, we have
$$
   \bigg\Vert f-\sum_{k\in \Z}2^j\int\limits_{[-2^{-j-1},2^{-j-1}]}f(-2^{-j}k+t)\,dt\,\phi(2^j\cdot+k)\bigg\Vert_p\le  C \omega_4(f,2^{-j})_p.
$$

\bigskip

\textbf{Example 2.} Now let $d=2$ and $U=[-1/2,1/2]^2$. Simple calculations show that in this case one has
$$
T(\xi)=g_0(\xi_1)g_0(\xi_2)+\frac{\pi^2}3g_2(\xi_1)g_0(\xi_2)+\frac{\pi^2}3g_0(\xi_1)g_2(\xi_2)=1+\frac{\pi^2}3(g_2(\xi_1)+g_2(\xi_2)),
$$
and, therefore,
\begin{equation*}
  \begin{split}
     \vp(x_1,x_2)=\frac56\sinc x_1\sinc x_2&+\frac{\sinc x_2}{24}(5\sinc(x_1+1)-4\sinc(x_1+2)+\sinc(x_1+3))\\
     &+\frac{\sinc x_1}{24}(5\sinc(x_2+1)-4\sinc(x_2+2)+\sinc(x_2+3)).
   \end{split}
\end{equation*}
It follows  from~\eqref{eqcor1++} that
\begin{equation*}
   \bigg\Vert f-\sum_{k\in\Z^2}m^j\int\limits_{M^{-j}[-1/2,1/2]^2}f(-M^{-j}k+t)\,dt\,\phi(M^j\cdot+k)\bigg\Vert_p\le  C \omega_4(f,\|M^{-j}\|)_p.
\end{equation*}

\bigskip

\textbf{Example 3.} Similarly to Example 2, if $d=2$ and  $U=B_1$,  taking into account that
$$
\w\vp (x)=\frac{\Gamma(1+d/2)}{\pi^{d/2}}\chi_{B_1}(x),\quad  \h{\widetilde{\vp}}(\xi)=\Gamma(1+d/2)\frac{J_{d/2}(2\pi|\xi|)}{(\pi|\xi|)^{d/2}},
$$
where $J_\lambda$  is the Bessel function of the first kind of order $\lambda$,
we obtain
$$
T(\xi_1,\xi_2)=1-\pi^2(g_2(\xi_1)+g(\xi_2))
$$
and
\begin{equation*}
  \begin{split}
     \vp(x_1,x_2)=\frac32\sinc x_1\sinc x_2&-\frac{\sinc x_2}{8}(5\sinc(x_1+1)-4\sinc(x_1+2)+\sinc(x_1+3))\\
     &-\frac{\sinc x_1}{8}(5\sinc(x_2+1)-4\sinc(x_2+2)+\sinc(x_2+3)).
   \end{split}
\end{equation*}
Hence,
\begin{equation*}
   \bigg\Vert f-\sum_{k\in\Z^2}\frac{m^j}{\pi}\int\limits_{M^{-j} B_1}f(-M^{-j}k+t)\,dt\,\phi(M^j\cdot+k)\bigg\Vert_p\le  C \omega_4(f,\|M^{-j}\|)_p.
\end{equation*}


\bigskip

III. Another improvement of the estimate given in Proposition~\ref{pr1 } can be  obtained by using an appropriate linear combination of the averaging operator rather than linear combinations of the function $\vp(x)=\sinc(x)$. Thus, the following estimate is a trivial consequence of Proposition~\ref{pr2}:

\begin{equation}\label{103}
  \bigg\Vert f-\sum_{k\in\zd}\sum_la_l\frac{m^j}{\mes U}\int\limits_{M^{-j}(U+l)}f(-M^{-j}k+t)\,dt\,\sinc(M^j\cdot+k)\bigg\Vert_p\le  C \omega_n(f,\|M^{-j}\|)_p,
\end{equation}
where $a_l$ is the $l$-th coefficient of the polynomial $T$ defined by~(\ref{102}) and the constant $C$ does not depend on $f$ and $j$.

Now let us obtain an analog of~(\ref{103}) for other functions $\phi$.

\begin{prop}
\label{pr3}
Let $f\in L_p$, $1<p<\infty$, $n\in \N$, and let  $U$ be a bounded measured subset in $\rd$. Suppose the function $\phi$ satisfies conditions of Theorem~\ref{th2}. Then there exists a finite set of numbers $\{b_l\}_{l\in \Z^d}\subset \C$ depending only on $d$, $n$, $U$, and $\vp$ such that
\begin{equation}\label{107}
  \bigg\Vert f-\sum_{k\in\zd}\sum_lb_l\frac{m^{j}}{\mes U}\int\limits_{M^{-j}(U-l)}f(-M^{-j}k+t)\,dt\,\phi(M^j\cdot+k)\bigg\Vert_p\le  C \omega_n(f,\|M^{-j}\|)_p,
\end{equation}
where the constant $C$ does not depend on $f$ and $j$.
\end{prop}

\textbf{Proof.}
Find complex numbers $c'_\alpha$, $\alpha\in\zd_+$, $[\alpha]< n$, satisfying
$$
c'_\nul=1, \quad\sum\limits_{\nul\le\alpha\le\beta}\left({\beta\atop\alpha}\right)
\overline{D^{\beta-\alpha}\h{\phi}(\nul)}c'_\alpha=0\quad \forall \beta\in\zd_+,\, \nul<[\beta]< n.
$$
Next, for $\w\vp (x)=\frac1{\mes U}\chi_U(x)$, we find  complex numbers $c_\alpha$, $\alpha\in\zd$, $[\alpha]< n$, satisfying
$$
c_\nul=1, \quad\sum\limits_{\nul\le\alpha\le\beta}\left({\beta\atop\alpha}\right)
\overline{D^{\beta-\alpha}\h{\w\phi}(\nul)}c_\alpha=c'_\beta\quad \forall \beta\in\zd_+,\, \nul<[\beta]< n,
$$
and set
\be
Q(\xi)=\sum_{\nul\le[\alpha]\le n}c_\alpha \prod_{j=1}^dg_{\alpha_j}(\xi_j),
\label{108}
\ee
where $g_\alpha$ is as in~(\ref{102}). If now $Q(\xi)=\sum\limits_l b_l e^{2\pi i(l,\xi)}$, then setting
$$
  \w\psi(x)=\sum_l b_l \w\phi(x+l),
$$
we obtain that  $D^{\beta}(1-\overline{\h\phi}{\h{\w\psi}})(\textbf{0}) = 0$ for all $\beta\in\zd_+$, $[\beta]<n$. Thus, due to Theorem~\ref{th2}, we have~\eqref{107}.~$\Diamond$

\bigskip

\textbf{Example 4.} In this example, we consider radial functions $\vp(x)=R_\delta(x)$ given by the
Bochner-Riesz type kernel
$$
R_\delta(x):=\frac{\Gamma(1+\delta)}{\pi^\delta}\frac{J_{d/2+\delta}(2\pi|x|)}{|x|^{d/2+\delta}}.
$$
Some results on approximation properties of sampling expansions generated by radial functions with a diagonal matrix $M$ can be found in~\cite{BFS}.
Let also $d=2$, $n=4$, and $U=B_1$.
By analogy with the above examples, we can compute that the polynomial $Q(\xi)$ given in~\eqref{108} has the form
$$
Q(\xi)=1+(2\delta-\pi^2)(g_2(\xi_1)+g_2(\xi_2))
$$
and, therefore, by~\eqref{107}, we derive
\begin{equation*}
  \bigg\Vert f-\sum_{k\in \Z^2}\sum_{l_1=0}^3 \sum_{l_2=0}^3 \frac{b_{l_1,l_2} m^{j}}{\pi}\!\!\!\!\!\!\!\int\limits_{M^{-j}(B_1-(l_1,l_2))}f(-M^{-j}k+t)\,dt\,R_\delta(M^j\cdot+k)\bigg\Vert_p\le  C \omega_4(f,\|M^{-j}\|)_p,
\end{equation*}
where
$$
b_{0,0}=1-\frac{2\delta-\pi^2}{2\pi^2},\quad
b_{1,0}=b_{0,1}=\frac{5(2\delta-\pi^2)}{8\pi^2},
$$
$$
b_{2,0}=b_{0,2}=\frac{-(2\delta-\pi^2)}{2\pi^2},\quad
b_{3,0}=b_{0,3}=\frac{2\delta-\pi^2}{8\pi^2},
$$
and
$$
b_{1,1}=b_{1,2}=b_{2,1}=0.
$$

\bigskip

IV.
Finally, from Theorem~\ref{thJack}, we obtain the following result related to the classical Kotelnikov decomposition.
\begin{prop}
Let $f\in L_p$, $1<p<\infty$,  $\h f$ be locally summable, and $n\in \N$. Then
\begin{equation}\label{eqcor1+}
  \bigg\Vert f- m^{j/2}\sum_{k\in \Z^d}\int\limits_{[-1/2,1/2]^d}  \h f({M^*}^j\xi)e^{-2\pi i (k,\xi)}d\xi\,\sinc(M^j\cdot+k) \bigg\Vert_p\le C\omega_n\(f,\|M^{-j}\|\)_p,
\end{equation}
where the constant $C$ does not depend on $f$ and $j$.
\end{prop}

\textbf{Proof.}
We apply Theorem~\ref{thJack} for $\vp(x)=\w\vp(x)=\sinc(x)$.
Since  $\h {\sinc}(\xi)=\chi_{[-1/2,1/2]^d}(\xi)$ and
$
\h {\sinc_{jk}}(\xi)=m^{-j/2}e^{2\pi i (k,{M^*}^{-j}\xi)}\chi_{[-1/2, 1/2]^d}({M^*}^{-j}\xi),
$
we have
$$
\langle f,\sinc_{jk}\rangle=\langle \h f,\h {\sinc_{jk}}\rangle=m^{j/2}\int\limits_{\rd}\h f({M^*}^j\xi)\overline{ \h{\w\phi}({M^*}^j\xi)}d\xi=\int\limits_{[-1/2,1/2]^d}  \h f({M^*}^j\xi)e^{-2\pi i (k,\xi)}d\xi,
$$
which,  by Theorem~\ref{thJack}, proves the proposition.
$\Diamond$

\medskip

Finally, let us note that~\eqref{eqcor1+} can be also written in the following form
\begin{equation}\label{eqcor2+}
  \bigg\Vert f-\sum_{k\in \Z^d}\mathcal{F}^{-1}(\chi_{M^{*j}[-1/2,1/2]^d}\h f)(M^{*-j}k)\,\sinc(M^j\cdot+k) \bigg\Vert_p\le C\omega_n\(f,\|M^{-j}\|\)_p.
\end{equation}

%


\medskip





\begin{thebibliography}{99}

\bibitem{BDR} {\sc C. de Boor, R. DeVore, A. Ron}, Approximation from shift-invariant subspaces of $L_2(\R^d)$, {\it Trans. Amer. Math. Soc.} \textbf{341} (1994), no.~2, 787--806.

\bibitem{DB-DV1} {\sc C. de Boor, R. DeVore, A. Ron},
Approximation orders of FSI spaces in $L_2(\rd)$,
{\it Constr. Approx.} \textbf{14} (1998), no.~3, 411--427.


%
%












\bibitem{BBSV} {\sc C. Bardaro, P.\,L. Butzer, R.\,L. Stens, G. Vinti}, Kantorovich-type generalized sampling series in the setting of Orlicz spaces, \emph{Sampl. Theor. Signal Image
Process.} \textbf{6} (2007), 29--52.

\bibitem{BM} {\sc C. Bardaro, I. Mantellini}, On pointwise approximation properties of multivariate semi-discrete sampling type operators, \emph{Results in Math.} (2017), DOI 10.1007/s00025-017-0667-7.

\bibitem{Brown}
        {\sc J.\,L. Brown, Jr.},
        On the error in reconstructing a non-bandlimited function by means of
        the bandpass sampling theorem,
        {\it J. Math. Anal. Appl.} {\bf 18} (1967), 75--84.





\bibitem{BFS}
{\sc P.\,L. Butzer, A. Fisher, R.\,L. Stens}, Generalized sampling approximation
of multivariate signals: theory and applications, {\it Note di Matematica} {\bf 10} (1990), no. 1,
173--191.








\bibitem{CS} {\sc D. Costarelli, R. Spigler}, Convergence of a family of neural network operators of the Kantorovich type, \emph{J. Approx.
Theory} \textbf{185} (2014), 80--90.



\bibitem{CV2} {\sc D. Costarelli, G. Vinti}, Rate of approximation for multivariate sampling Kantorovich operators
on some functions spaces, {\it J. Int.
Eq. Appl.} \textbf{26} (2014), no.~4,  455--481.


\bibitem{CV0} {\sc D. Costarelli, G. Vinti},
Degree of approximation for nonlinear multivariate sampling Kantorovich operators on some functions spaces,
\emph{Numer. Funct. Anal. Optim.} \textbf{36} (2015), no.~8, 964--990.




\bibitem{FCMV} {\sc F. Cluni, D. Costarelli, A.\,M. Minotti, G. Vinti},
Applications of sampling Kantorovich operators to thermographic images for seismic engineering,
\emph{J. Comput. Anal. Appl. }\textbf{19} (2015), no.~4, 602--617.





%



%

 \bibitem{v58}
{\sc R.-Q. Jia},  Refinable shift-invariant spaces: From splines to wavelets.
[CA] Chui, C. K. (ed.) et al., Approximation theory VIII. Vol. 2.
Wavelets and multilevel approximation. Papers from the 8th Texas international
conference, College Station, TX, USA, January 8--12, 1995. Singapore:
World Scientific. Ser. Approx. Decompos. 1995, Vol. 6, pp.~179--208.


\bibitem{Jia1} {\sc R.-Q. Jia}, Convergence rates of cascade algorithms, {\it Proc. Amer. Math. Soc.} \textbf{131} (2003), 1739--1749.

\bibitem{Jia2}
{\sc R.-Q. Jia}, Approximation by quasi-projection operators in Besov spaces,
{\it J. Approx. Theory} \textbf{162} (2010), no.~1, 186--200.




\bibitem{KM} {\sc A. Kivinukk, T. Metsmagi}, On boundedness inequalities of some
semi-discrete operators in connections with sampling operators,
Proceedings of 2015 International conference on Sampling Theory
and Applications, Washington DC, USA, 25-29 May 2015, IEEE
Xplore, pp.~48--52.


\bibitem{KKS}
{\sc Yu. Kolomoitsev, A. Krivoshein, M. Skopina}, Differential and falsified sampling  expansions, arXiv:1703.10420v1  [math.CA]  30 Mar 2017.


\bibitem{KS}
{\sc A. Krivoshein, M. Skopina}, Approximation by frame-like wavelet systems,
{\it Appl. Comput. Harmon. Anal.}
\textbf{31} (2011), no.~3, 410--428.


\bibitem{KS1}
{\sc A. Krivoshein, M. Skopina}, Multivariate sampling-type approximation,
{\it Anal. Appl.} (2016),
DOI: 10.1142/S0219530516500147.


\bibitem{KPS}
{\sc A. Krivoshein, V. Protasov, M. Skopina}, Multivariate Wavelet Frames, Springer, 2016.

\bibitem{LJC}
{\sc J. J. Lei, R. Q. Jia, E. W. Cheney}, Approximation from
shift-invariant spaces by integral operators, {\it SIAM J. Math. Anal.}
\textbf{28} (1997), 481--498.

\bibitem{Nik} {\sc S. M. Nikol'skii},
The Approximation of Functions of Several Variables and the Imbedding Theorems, 2nd ed.
{Moscow}: Nauka, 1977 (Russian). -- English transl. of 1st. ed.: John Wiley $\&$ Sons, New-York, 1978.

\bibitem{NPS} {\sc  I. Ya. Novikov, V. Yu. Protasov,  M. A. Skopina}, Wavelet Theory, Amer. Math. Soc, 2011.



\bibitem{OT15} {\sc O. Orlova, G. Tamberg},  On approximation properties of generalized
Kantorovich-type sampling operators, \emph{J. Approx. Theory} {\bf 201} (2016), 73--86.




%
				

\bibitem{Sk3}				
{\sc M. Skopina}, On construction of multivariate wavelets with vanishing moments, \emph{Appl. Comput. Harmon. Anal.} \textbf{20} (2006), no.~3, 375--390.


\bibitem{Sk1}
         {\sc M. Skopina},
        Band-limited scaling and wavelet expansions,
        {\it Appl. Comput. Harmon. Anal.} {\bf 36} (2014), 143--157.






\bibitem{Timan} {\sc A. F. Timan}, Theory of Approximation of Functions of a Real Variable, Pergamon Press, Oxford, London, New York, Paris, 1963.


\bibitem{VZ1} {\sc G. Vinti, L. Zampogni}, Approximation by means of nonlinear Kantorovich sampling type operators in Orlicz spaces,
\emph{J. Approx. Theory} \textbf{161} (2009), 511--528.

\bibitem{VZ2} {\sc G. Vinti, L. Zampogni}, Approximation results for a general class of Kantorovich type operators, \emph{Adv. Nonlinear
Stud.} \textbf{14 }(2014), no.~4, 991--1011.


\bibitem{Vladimirov-1}
 {\sc V.~S. Vladimirov},
Generalized Functions in Mathematical Physics,
MIR, 1979 (translated from Russian).


\bibitem{Wil} {\sc  G. Wilmes},  On Riesz-type inequalities and $K$-functionals related to Riesz potentials in $\R^N$,
{\it Numer. Funct. Anal. Optim.} {\bf 1} (1979), no.~1, 57--77.

\end{thebibliography}
\end{document}